\documentclass[11pt]{article}
\usepackage{graphicx}
\newtheorem{claim}{Claim}

\title{A simple introduction to Karmarkar's Algorithm 
for Linear Programming}

\author{Sanjeev Saxena\thanks{E-mail: ssax@iitk.ac.in}\\
Dept. of Computer Science and
Engineering,\\ Indian Institute of Technology,\\
Kanpur, INDIA-208 016}

\date{December 21, 2017}

\begin{document}
\maketitle

\subsection*{\centering{Abstract}}

An extremely simple, description of Karmarkar's algorithm with very
few technical terms is
given.

\section{Introduction}

A simple description of Karmarkar's algorithm\cite{karm} together with
analysis is given in this paper. Only knowledge of simple algebra, vector
dot product and matrices is assumed.  Even though the method is described
in several books \cite{taha,bazar,hand,Kar,Sai}, analysis is either left
out \cite{taha} or is fairly complicated. We show that 
the presentation of Roos,  Terlaky and Vial\cite{RT} can be further simplified,
(mainly) by using only essential notation. Allmost all
results in this paper are from
standard textbooks (see References).

Let $A$ be an $m\times n$ matrix of rank $m$ and $e$ is a vector of all
ones. Karmarkar's Problem is:\\
$\min c^Tx$ subject to constraints 
$Ax=0$, $e^Tx=1$ and $x\ge 0$.

Further, it is assumed that optimal value of $c^Tx$ is zero and all ones
vector is feasible, i.e., $Ae=0$.
We have to either find a point of cost $0$ or show that none exist

We first scale the variables: $x'=xn$, then the problem becomes $\min
c^T(x'/n)$ or equivalently $\min c^T(x')$ subject to $Ax'=0$,
$e^Tx'=n$ and $x'\ge 0$. We will drop the primes, and the problem
is\cite{RT}:

$\min c^Tx$ subject to constraints
$Ax=0$, $e^Tx=n$ and $x\ge 0$.

Remark: Problem is trivial if $c^Te=0$, hence we assume that $c^Te>
0$.

We will assume that all $m+1$ equality constraints are linearly
independent (else we eliminate redundant rows of $A$).

We need few definitions.
{\textbf{Standard simplex}} consists of points in $n$ dimensions s.t.
$e^Tx=n$, $x\geq 0$.
The centre of the simplex is $e=(1,1,{\ldots} ,1)$.

Let $R$ be the radius of {\textbf{outer sphere}}, the smallest sphere
containing the standard simplex (circumscribes standard simplex). 
$R$ is the distance from $e$ to one of the corner point (see Figure for
an example in three dimensions) of the standard simplex (say,
$(0,0,{\ldots} ,n)$), thus $R^2=(0-1)^2+{\ldots} +(0-1)^2+(n-1)^2=
(n-1)^2+ (n-1)=n(n-1)$, or $R=\sqrt{n(n-1)}$.

\begin{figure}
\includegraphics[height=1.8in]{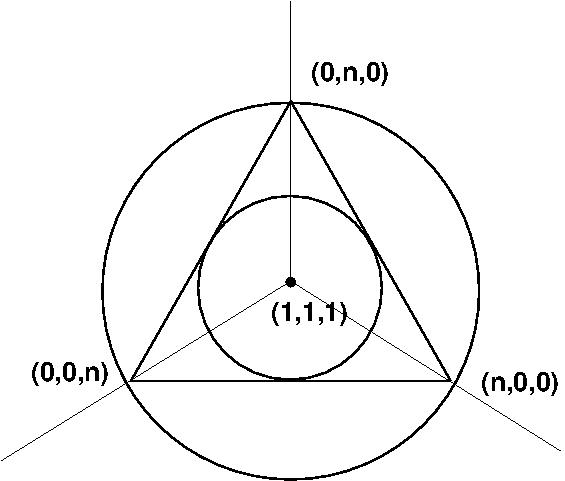}
\end{figure}

Let $r$ be the radius of another sphere--- the largest sphere centred
at $e$ and completely inside the standard simplex (inscribed inside the
standard simplex). 
This sphere will be tangent to each face. Each face will have one
coordinate as $0$. By symmetry, all other coordinates at the point of
contact will be same (say $w$). As the point of contact is on the
standard simplex $S$, $0+w+w+{\ldots} +w=n$ or $(n-1)w=n$ or $w=n/(n-1)$.
Hence,
$$r^2=(0-1)^2+(w-1)^2+{\ldots} +(w-1)^2=1+(n-1)(w-1)^2=1+\frac{1}{n-1}
=\frac{n}{n-1}$$
Or $$r=\sqrt{\frac{n}{n-1}}$$

We take $e=(1,1,{\ldots} ,1)$ as the starting point (which by assumption
is feasible). Then, we minimise the objective function over a smaller
sphere, which we will call the {\textbf{inner sphere}}, having same centre
$e$, but radius $\alpha r$, 
less than $r$ (we will see in
Section~5.2 that $\alpha$ can be chosen as $1/(r+1)$). Let us assume that
the minimum occurs at point $z$. Then, in next iteration, we take the
starting 
point 
as $z$.  The
problem of minimisation on sphere is discussed in next section
(Section~2). Finally, the point $z$ is mapped to $e$ and the process
repeated; the details of mapping are in Section~4.
In Section~3, it is shown that the objective value at next point $z$ is a
fraction of that at the initial point $e$.

\section{Mathematical Preliminaries}

Let $P$ be a point in (possible hyper) plane and $\hat{n}$ a unit vector
normal to it. If $P_0$ is any point in plane, then the (vector) dot
product will be zero, i.e., $PP_0\bullet \hat{n}=0$. Writing in full,\\
$n_1(x_1-x'_1)+n_2(x_2-x'_2)+{\ldots} +n_k(x_k-x'_k)=0$
or equivalently,\\
$n_1x_1+n_2x_2+{\ldots} +n_kx_k=n_1x'_1+n_2x'_2+{\ldots} +n_kx'_k=C$
Thus, for $a\bullet x=b$, normal will be in the direction of vector $a$.

Equation of sphere with centre $Q$ as $\beta$ is: $(x-\beta)\bullet
(x-\beta)=r^2$
Let $P$ be any point in plane $a\bullet x=b$, then we know that $a\bullet
x_P=b$. If we want $QP$ to be perpendicular to plane, then
$(x_P-\beta)=\mbox{(constant)} a$

\begin{figure}
\includegraphics[height=1.8in]{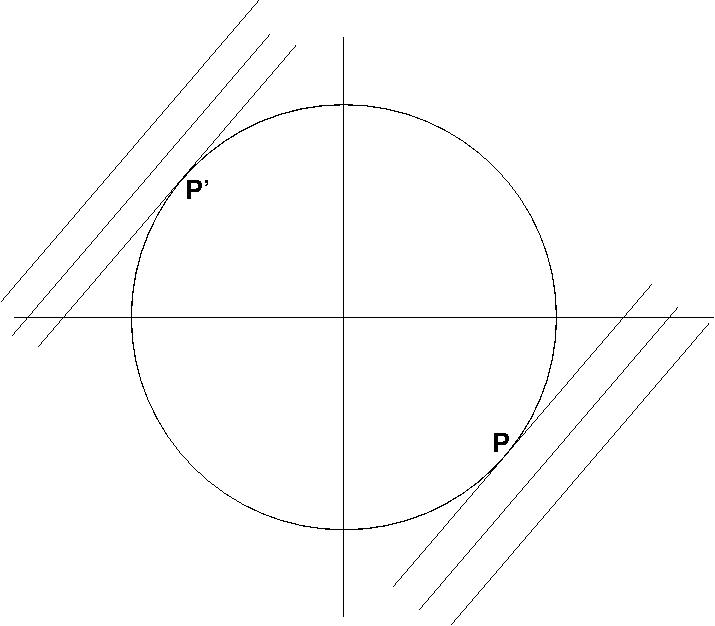}
\end{figure}

Consider the problem \cite{Murty}:
$\min c^Tx$ subject to $\sum (x-a)^2\leq \rho^2$.
There are two points on the sphere (see figure for an example in two
dimensions), where planes parallel to given plane can be tangent--- one
corresponding to maxima and other minima. These are $+n$ and $-n$. These
points (say ``$P$'') are in direction $\pm \hat{n}$ (a unit vector in
direction of $c$) with length $|x_P-\beta|=r$.
Thus, to minimise $c^T x$ over our sphere, we start at the centre $\beta$
and take a step of length $r$ in direction $-c^T$
Less informally, if $c=0$, all points on the sphere are optimal. If
$c\neq 0$, the solution is obtained by taking a step of length $\rho$
(radius of sphere) from the centre $a$ in the direction $-c$; this can be
seen by considering ``parallel'' planes $c^Tx=$constant. The point
at which minimum will be attained will be the point of contact on the
sphere to a tangent plane.

Next \cite{Murty} consider the problem: $\min c^Tx$ subject to $Ax=b$ and 
$\sum (x-a)^2\leq \rho^2$.
If $c=0$, all points common to (i.e., on the intersection of) sphere and
the plane are optimal. If $c\neq 0$, let $\bar{c}$ be the orthogonal
projection of $c$ onto the plane $Ax=b$. If $\bar{c}=0$, then $c$ is
a linear combination of rows of $A$ and the objective function (on the
intersection of sphere and plane) is constant, and all ``feasible''
points are also optimal. If $\bar{c}\neq 0$, then the solution is
obtained by taking a step of length equal to the radius of
lower-dimensional sphere (intersection of our sphere and plane
$Ax=b$) from the centre $a$ in the direction $-\bar{c}$.

\section{Solution on inner Sphere}

Let us consider the problem: $\min c^Tx$ subject to $Ax=0$, $\sum_i
x_i=n$ and $\sum_i (x_i-1)^2< \alpha^2 r^2$, i.e, we are minimising only
over points of the inner sphere. As the inner sphere is completely inside
the standard simplex, we can drop the constraint $e^Tx=n$ (and also $x
\geq 0$), and the problem becomes\\
$\min c^Tx$ s.t. $Ax=0$ and $\sum
(x_i-1)^2< \alpha^2r^2$. 

As the point $e$ is feasible, $Ae=0$, $e$ lies on the plane $Ax=0$.
Moreover, as $e$, the centre of inner sphere, 
lies on $Ax=0$, the intersection of $Ax=0$ and the inner sphere will be a
``sphere'' of the same radius $\alpha r$ but in lower dimension. (e.g.,
intersection of a sphere and a plane containing the centre is a circle
with same radius and centre).

The minimum value of a linear function on sphere will be at a point
(lower one) where linear function touches the sphere. The radius vector
at that point will be perpendicular to the plane. Let us assume that
$\hat{p}$ is a unit vector in that direction. Then the point at which
minimum will be attained will be $e-\hat{p}(\mbox{radius})$. Thus, for
outer sphere it will be (say) $z_R=e-\hat{p}R$ and for the inner sphere
it will
be 
$z_\alpha=e-\hat{p}\alpha r$. We will see how to determine $\hat{p}$
later (see Section~4.2).

As the outer sphere completely contains the solution space, the minimum
value (of objective function $c^Tx$) will be smaller than (or equal to)
the actual optimal value which is zero. Thus, (assuming minimum value is
at $z_R$)

$$0\geq c^Tz_R=c^T(e-\hat{p}R) \mbox{ or } c^T\hat{p}R \geq c^Te$$

Thus,  
$c^T\hat{p}\geq \frac{c^Te}{R}$

As the inner sphere is inside the simplex, value of the objective
function can not be less than that on the simplex. But, as the optimal
value on simplex is zero, optimal value on the inner sphere is
non-negative.
If the minimum value occurs at $z_\alpha$ on the inner sphere, then 
the value of objective function 
is: 

$$0\leq c^Tz_\alpha=c^T(e-\hat{p}\alpha r) \mbox{ or } c^T\hat{p}\alpha r \leq c^Te$$
But as $\frac{c^Te}{R}\leq c^T\hat{p}$, we get
$$c^Tz_\alpha=c^Te-c^T\hat{p} \alpha r\leq
c^Te-\left(\frac{c^Te}{R}\right)\alpha r=c^Te\left(1-\frac{\alpha
r}{R}\right)$$

Thus, if we start with initial solution $e$, then $z_\alpha$ (the next
solution) is an improvement (in value of objective function) by a factor
of $\left(1-\frac{\alpha r}{R}\right)$ over the initial solution.

\section{Karmarkar Transform and Algorithm}

We will like to map the new point $z$ again to $e$ to repeat the process.
Thus, for 
$a=(a_1,a_2,{\ldots} ,a_n)$ and $x=(x_1,x_2,{\ldots} ,x_n)$, we define
the transform $y=T_a(x)$ with (for $i=1,{\ldots} ,n$)

$$y_i=n \frac{a_ix_i}{\sum_j a_jx_j}$$

If all $a_i$s and $x_i$s are non-negative (at least $a_ix_i$ should be
non-zero), each component of transform will be less than $n$ and sum of
all coordinates will be $n$, thus the range is again our standard simplex. 
Moreover, if $b=a^{-1}=\left(\frac{1}{a_1},\frac{1}{a_2},{\ldots} ,
\frac{1}{a_n}\right)$, then $T_b$ will be the inverse transformation.
Thus this transformation is one-to-one on our standard simplex. 

Remark: This can also be seen directly. If $T_a(x)=T_a(y)$, then
$\frac{x_1}{y_1}= \frac{x_1}{y_2}={\ldots} =\frac{x_n}{y_n}= r$ (say).
But as $\sum x_i=\sum y_i=1$, we have $r=1$.

Moreover, if $\lambda$ is any number then if $y'=T_a(\lambda x)$ then 
$$y'_i = \frac{n \lambda a_ix_i}{\sum \lambda a_jx_j} 
= \frac{n {\lambda} a_ix_i}{{\lambda} \sum a_jx_j}
= \frac{n a_ix_i}{\sum a_jx_j}\\
= y_i$$
Thus, $T_a(\lambda x)=T_a(x)$.

\subsection{Modified Problem}

Assume that we are applying transform $T_a$ (with $a_i=1/z_i$), then point
$z$ will be mapped to $e$. Let $D$ be a diagonal matrix with diagonal
entries: $D=\mbox{ diag}(a_1,{\ldots} ,a_n)$. 

If $z$ is any feasible point, then for any positive $x> 0$ satisfying
$e^Tx=1$, we saw that there is a unique point $\xi$ (which depends on
$x$) s.t. $x=T_z(\xi)$. 

Remark: $x=T_z(\xi)$ implies $x_i=z_i\xi_i/(\sum z_i\xi_i)$,
thus $\xi_i=\rho \frac{x_i}{z_i}$ where $\rho$ is such that $\sum
\xi_i=1$.

Then, equation $Ax=0$ is equivalently to $\sum_j A_{ij}x_j=0$ or $\sum_j
A_{ij}\rho x_j=0$ or equivalently, $\sum A_{ij} \xi_j z_j=0$ or $A
\xi z=0$.

The objective function $c^Tx=\sum c_ix_i=\frac{1}{\rho}\sum c_i\xi_i
z_i$. As optimal value of $c^Tx$ is zero, it follows that the optimal
value of the transformed objective function $\sum c_i\xi_i z_i$ is also
$0$. 

Replacing $\xi$ by $x'$ and using $Z=\mbox{diag}(z_1,{\ldots} ,z_n)$, 
the transformed problem is:
\begin{quote}
$\min (Zc)^Tx'$ subject to $AZx'=0$, $e^Tx'=n$ and $x'\geq 0$.
\end{quote}
Moreover, as $z$ is a feasible point $Az=0$ or equivalently $AZe=0$, thus
$e$ is again feasible.

We can thus repeat the previous method with $Zc$ instead of $c$ and $AZ$
instead of $A$. In other words, we have to minimise the modified
objective function over inscribed sphere of radius  $\alpha
r$ i.e., the ``modified'' problem is:
\begin{quote}
$\min (Zc)^Tx$ subject to $AZx=0$, $e^Tx=n$ and $||x-e||\le \alpha r$.
\end{quote}

\subsection{Result from Algebra-- Projection Matrix} 

Assume that $A$ is $m\times n$ matrix, then rank of $A$ is said to be
$m$, with $m<n$ iff all $m$ rows of $A$ are linearly independent, i.e.,
$\beta_1A_1+\beta_2A_2+{\ldots} +\beta_mA_n=0$ (here $0$ is a row
vector of size $n$) has only one solution $\beta_i=0$. Thus, if $v$ is
any $1\times m$ matrix (a column vector of size $m$), then $vA=0$ implies
$v=0$.  We use the fact that the matrix $AA^T$ has rank $m$ and is invertible\footnote{{
As $A$ is $m\times n$ matrix, $A^T$ will be $n\times m$ matrix. The
product $AA^T$ will be an $m\times m$ square matrix. Let $y^T$ be an
$m\times 1$ matrix (or $y$ is a row-vector of size $m$). 

Consider the equation $(AA^T)y^T=0$. Pre-multiplying by $y$ we get
$yAA^Ty^T=0$ or $(yA)(yA)^T=0$ or the dot product $<yA,yA>=0$ which, for
real vectors (matrices) means, that each term of $yA$ is (individually)
zero, or $yA=0$, which implies $y$ is identically zero.
Thus, the matrix $AA^T$ has rank $m$ and is invertible.}}.

Let $A$ be an $m\times n$ matrix, then $Ax=b$ (hence $Ax=0$) represents a
set of $m$ equations. Each equation will be a hyperplane in $n$
dimensions. Let $v$ be a vector of size $n$, we will like to ``project''
$v$  onto the (lower dimension or intersection of) hyperplane $Ax=b$. If
$P=pv$ is the projection, then \cite{algebra}, we wish to write $p=Pv$,
as best as possible, as $p=Pv=\alpha_1a_1+\alpha_2a_2+{\ldots}
+\alpha_ma_m$ where $a_1,{\ldots} ,a_m$ are rows of $A$. Projection $p$
will again be a vector of size $n$. This can be written as 
$p=A^T\alpha$.

The error of projection (or rejection) $E=v-p=v-Pv$ will also be a vector
of size $n$.  Then as error $E$ can not be in these hyperplane, i.e., it
should be ``perpendicular'' to these hyperplanes, thus we want $A E=0$. Or
$A(v-p)=0$, or $A(v-A^T\alpha)=0$, thus $Av=AA^T \alpha$, or
$\alpha=(AA^T)^{-1}Av$. 

Hence, $p=Pv=A^T\alpha$, or $Pv=A^T(AA^T)^{-1}Av$, or $P=A^T(AA^T)^{-1}A$.
Or, $E=v-Pv=v-A^T(AA^T)^{-1}Av=\left(I-A^T(AA^T)^{-1}A\right)v$

Thus, to summarise,
projection matrix is 
$P=A^T(AA^T)^{-1}A$ and
Rejection matrix (to get the part perpendicular to the hyperplanes) is:
$I-A^T(AA^T)^{-1}A$

\subsection{Algorithm}

Start the algorithm with $x=e$ and find $z$, the optimal value on the
inscribed sphere of radius $\alpha r$ (we will see more details of this
step later).

In a typical iteration, we next apply a transformation $T_{z^{-1}}$ to
map $z$ to $e$ (we need a transformation $T_b$ such that $z$ gets mapped to
$e$, then $b=1/z$). Modify $A$ and $c$ and find a new $z$ (say $z'$), the
optimal value on the inscribed sphere of radius $\alpha r$.

Remark: As we know that transformed value of $z$ is $e$, transform
$T_{z^{-1}}$ is not actually applied, only the values of $A$ and $c$ are
updated.

Then apply inverse transformation $T^{-1}_{z^{-1}}$ (or $T_z$) to map
$z'$ to original space, to get the value of $x'$ for next iteration.

The formal algorithm is (value of $\alpha$ is fixed to $1/(r+1)$ in
Section~5.2):

Initialise: $r=\sqrt{\frac{n}{n-1}}$, $\alpha=\frac{1}{r+1}$,
$x=e=(1,{\ldots} ,1)$

Main step: If $c^Tx<\epsilon$, return current $x$ as solution of desired
accuracy.

Let $D=\mbox{diag}(x_1,{\ldots} ,x_n)$

\begin{displaymath}
P=\left(
\begin{array}{c}
AD\\
I
\end{array} \right)
\end{displaymath}

Let 
\begin{eqnarray*}
c_P &=& \left(1-P^T\left(PP^T\right)^{-1}P\right)(cD)^T\\
\hat{p} &=& \frac{c_P}{||c_P||}\\
z &=& e-\alpha r \hat{p}\\
x &=& \frac{nDy}{e^TDy}\\
\end{eqnarray*}

\section{Analysis---Potential Function}

For analysis, we define a potential function:
$$\Phi(x)= n \log c^Tx- \sum_{i=1}^n \log x_i$$

As $e^Tx=n$ we have $\frac{1}{n}x_i=1$. But as Arithmetic mean is greater
than or equal to geometric mean:

$1=\sum \frac{x_i}n\geq
\left(\prod x_i\right)^{1/n}$ or taking logs, $\frac{1}{n} \sum \log x_i
\leq 0$, or $\sum \log x_i \leq 0$, or $\Phi(x)\leq n \log c^Tx$ \mbox{
or } 
\begin{quote}
\fbox{$c^Tx \geq
\exp\left(\frac{\Phi(x)}{n}\right)$}
\end{quote}

Next observe that 
\begin{eqnarray*}
\Phi(\lambda x) &=& n \log c^T(\lambda x)- \sum_{i=1}^n
\log (\lambda x_i)\\
&=& \left(\log c^Tx+\log \lambda\right) -
\left(\sum_{i=1}^n \log x_i + \sum_{i=1}^n \log \lambda \right)\\
&=& n \log c^Tx- \sum_{i=1}^n \log x_i
=\Phi(x)
\end{eqnarray*}

Let $x$ be a positive vector in our standard simplex
and let $y=T_x(z)$, then
$$y_i=\frac{nx_iz_i}{\sum x_iz_i},\mbox{ or }\Phi(y)=\Phi(T_z(x))=\Phi(xz)$$

Moreover, $\Phi(xz)=n \log c^T(xz)-\sum_{i=1}^n \log x_iz_i=
n \log c^Txz- \sum \log x_i - \sum \log z_i$

Finally, 
\begin{eqnarray*}
\Delta \Phi &=&\Phi(x)-\Phi(y)
= \Phi(x)-\Phi(xz)\\
&=& n \log c^Tx- {\sum_{i=1}^n \log x_i}- 
  \left( n \log c^Txz- {\sum \log x_i} - \sum \log z_i \right)\\
&=& n\log \frac{c^Tx}{c^Txz}+\sum \log z_i
\end{eqnarray*}
But $\log \frac{c^Tx}{c^Txz}$ is the ratio of original and transformed
problems, which we saw reduces by at least
$\left(1-\alpha\frac{r}{R}\right)$. Thus:

$$\Delta \Phi \geq -n \log \left(1-\alpha\frac{r}{R}\right)+ \sum \log
z_i$$

We will show that $\Delta$ is more than a constant. But, before continuing
with analysis, we need some more results from algebra.

\subsection{Function $\Psi$}

In this section, all logs are to base $e$. Assume that $|x|<1$.
As
$\log (1+x)= x-\frac{x^2}{2}+\frac{x^3}{3}-\frac{x^4}{4}+{\ldots} $,
for $-1<y<0$, (let $y=-x$) we have\\
$y-\log(1+y)=-x-\log(1-x)=\frac{x^2}{2}+\frac{x^3}{3}+\frac{x^4}{4}+{\ldots}>0
$

We define\cite{RT} $\Psi(x)=x-\log (1+x)$.
Observe that $\Psi(0)=0$ 
Hence, $\Psi(x)=\frac{x^2}{2}-\frac{x^3}{3}+\frac{x^4}{4}-{\ldots} $ And,
$\Psi(x)\geq 0$ for $x>-1$

\begin{claim}
If $x>0$ and $a^2=b^2+c^2$, then if $ax,bx,cx$ are greater than $-1$ (for
$\log (1+hx)$ to be defined), 
then
$$\Psi(-|a|x)> \Psi(bx)+\Psi(cx)$$
\end{claim}

Proof: From $a^2=b^2+c^2$, we know that $|a|\geq |b|$ and $|a|\geq |c|$.
Without loss of generality, assume $a<0$, then $-|ax|=ax$ and we have

$\Psi(ax)= \frac{a^2x^2}{2}-\frac{a^3x^3}{3}+\frac{a^4x^4}{4}-{\ldots} $

$\Psi(bx)= \frac{b^2x^2}{2}-\frac{b^3x^3}{3}+\frac{b^4x^4}{4}-{\ldots} $

$\Psi(cx)= \frac{c^2x^2}{2}-\frac{c^3x^3}{3}+\frac{c^4x^4}{4}-{\ldots} $

As $a^2=b^2+c^2$, coefficients of $x^2$ on both sides are equal.

We thus need to compare $-\frac{a^3x^3}{3}$ with
$-\frac{(b^3+c^3)x^3}{3}$.
As $a^2=b^2+c^2$, $a^3=ab^2+ac^2$. Thus, we are comparing
$-x(ab^2+ac^2)$ and $-x(b^3+c^3)$ or equivalently
$0$ and $xb^2(a-b)+xc^2(a-c)$.

As $a<0$  and $-a>|b|$ and
$-a>|c|$, 
the expressions can be re-written as:
$0$ and $-xb^2(b-a)-xc^2(c-a)$.
Both terms inside brackets are positive, and as $x>0$, right hand side will be
negative. 

\begin{claim}
If $a^2=b^2+c^2+{\ldots} +s^2$, then if $ax,bx,cx, {\ldots} $ are greater
than $-1$ (for $\log (1+hx)$ to be defined),
then if  $x>0$, 
$\Psi(-|a|x) > \Psi(bx)+\Psi(cx)+{\ldots} +\Psi(sx)$
\end{claim}
Proof (induction): Let $v^2=c^2+{\ldots} +s^2$.  Then $a^2=b^2+v^2$.
Then from the above claim,
$\Psi(-|a|x) > \Psi(bx)+\Psi(-|v|x)$. 
By induction hypothesis,
$\Psi(-|v|x) > \Psi(cx)+{\ldots} +\Psi(sx)$
Proof follows by adding the two inequalities. 

\subsection{Analysis Continued}

As $-\log (1-x)= x+\frac{x^2}{2}+\frac{x^3}{3}+{\ldots} > x$, 
hence,
\begin{eqnarray*}
\Delta \Phi &\geq& -n \log \left(1-\alpha\frac{r}{R}\right)+ \sum \log
z_i\\
&=& n\alpha\frac{r}{R}+\sum \log z_i\\
&=& \alpha r^2+\sum \log z_i
\end{eqnarray*}
But $z=e-\alpha r \hat{p}$ or $z_i=1-\alpha r p_i$

And, $\log z_i=\log (1-\alpha r p_i)=-\left(-\alpha r p_i-\log (1-\alpha
r p_i)\right)-\alpha r p_i=-\alpha r p_i-\Psi(-\alpha r p_i)$.

But $\hat{p}$ is a unit vector, $\sum p^2_i=1$. As $z$ is feasible $\sum
z_i=n$, or $\sum p_i=0$. 
Hence,  $$-\sum \log z_i=\sum \alpha r p_i+\sum
\Psi(-\alpha r p_i)= \sum \Psi(-\alpha r p_i)\leq \Psi(-\alpha r)$$
The last inequality follows from  Corollary~1 above.

For $\Psi(-\alpha r)$ to be defined, $\alpha r<1$, we thus\footnote{{
Differentiating, $f(\alpha)=\alpha r^2-\Psi(-\alpha r)$ w.r.t. $\alpha$, 
$f'(\alpha)=
r^2-\frac{-\alpha r}{1-\alpha r}+n\frac{-r/R}{1-(\alpha r/R)}$
$f'(\alpha)=r^2+r+(-r)\frac{1}{1-\alpha r}$, equating $f'(\alpha)$ to
$0$, we get $1+r=\frac{1}{1-\alpha r}$, or $1-\alpha r=\frac{1}{1+r}$ or
$\alpha r=1-\frac{1}{1+r}=\frac{r}{r+1}$ or $\alpha=\frac{1}{1+r}$. Thus,
the maximum value is at $\alpha=\frac{1}{1+r}$. }} choose
$\alpha=\frac{1}{r+1}$, we have
\begin{eqnarray*}
\Delta \Phi &\geq& \alpha r^2-\Psi(-\alpha r)
= \alpha r^2-\big((-\alpha r)-\log(1-\alpha
r)\big)\\
&=&\alpha r^2+\alpha r+\log (1-\alpha r)\\
&=& \frac{r^2}{1+r}+\frac{r}{1+r}+\log \left(1-\frac{r}{r+1}\right)\\
&=& r+\log \left(1-\frac{r}{r+1}\right)\\
&=&r-\log (1+r)
= \Psi(r) \mbox{ but as } r=\sqrt{n/n-1}>1\\
&\geq& \Psi(1)
= 1-\ln 2\approx 0.3
\end{eqnarray*}

Hence, potential decreases by a fixed amount after each iteration.

After $k$, iterations, $\Phi(e)-\Phi(x)> k\Psi(1)$. But as $\Psi(e)=n
\log c^Te$, we get $\Phi(x)< n\log c^T e- k\Psi(1)$. As $x$ is inside the
standard simplex, $$c^Tx\leq \exp\left(\frac{\Phi(x)}{n}\right)<
\exp\left(\frac{n \log c^Te-k\Psi(1)}{n}\right)$$

If we stop as soon as error $c^Tx\leq \epsilon$, we get
\begin{eqnarray*}
\exp\left(\frac{n \log c^Te-k\Psi(1)}{n}\right) &\leq& \epsilon \mbox{ or}\\
\frac{n\log c^Te-k\Psi(1)}{n} &\leq& \log \epsilon \mbox{ or }\\
\log c^Te-k\Psi(1) &\leq&  n \log \epsilon \mbox{ or }\\
k &\geq& \frac{n}{\Psi(1)}\log \frac{c^Te}{\epsilon}
\end{eqnarray*}
Thus,  after at most
$\frac{n}{\Psi(1)}\log \frac{c^Te}{\epsilon}$ iterations, algorithm finds
a feasible point $x$ for which $c^Tx\le \epsilon$

\subsection*{Acknowledgments}

Many thanks to students of CS647A (2017-2018 batch) for their valuable
feedback, questions and comments on an earlier version.

\bibliography{general}

\end{document}